\documentclass[a4paper,11pt]{article}

\title{On Zermelo'-like problems: a Gauss-Bonnet inequality and a E. Hopf theorem}

\author{Ulysse Serres\footnote{Universit\`a degli Studi di Firenze, Dipartimento di Matematica per le Decisioni
via Lombroso 6/17, 50134 Firenze, Italy; email: \texttt{ulysse.serres@unifi.it}}}
\date{}

\usepackage{amsmath,amsthm,amssymb}
\usepackage{graphicx}
\usepackage{color}

\renewenvironment{proof}{\noindent{\emph{\textbf{P\hspace{-0.2mm}roo\hspace{-0.15mm}f.}}}}{\hfill$\blacksquare$}

\newtheorem{theorem}{Theorem}[section]
\newtheorem{definition}[theorem]{Definition}
\newtheorem{example}[theorem]{Example}
\newtheorem{lemma}[theorem]{Lemma}
\newtheorem{proposition}[theorem]{Proposition}
\newtheorem{coro}[theorem]{Corollary}
\newtheorem{remark}[theorem]{Remark}

\renewcommand{\le}{\leqslant}
\renewcommand{\ge}{\geqslant}
\renewcommand{\Vec}{{\rm Vec\,}}

\newcommand{\lang}{\left\langle}
\newcommand{\rang}{\right\rangle}
\newcommand{\bs}{\boldsymbol}
\newcommand{\parfrac}[2]{\frac{\partial{#1}}{\partial{#2}}}
\newcommand{\ad}{{\rm ad\, }}

\newcommand{\e}[2]{e^{#1\,#2}}
\newcommand{\vh}{\bs{\vec{h}}}

\newcommand{\he}{\bs{\hat{h}}}
\newcommand{\ve}{\bs{\hat{v}}}
\newcommand{\X}{\bs{X}}
\newcommand{\tg}{\tilde{g\,}\!}

\newcommand{\diffeo}[1]{F^{\sr#1\up}}
\newcommand{\sr}{\scriptscriptstyle}

\newcommand{\Z}{ {\sr \mathcal{Z}} }

\newcommand{\coZ}{ {\sr \rm{co}\mathcal{Z}} }
\newcommand{\phicoZ}{ {\varphi_g^{\sr \up}} }
\newcommand{\hZ}{h}
\newcommand{\vhZ}{\vh}

\newcommand{\hDZ}{h}
\newcommand{\vhDZ}{\vh}
\newcommand{\vhmag}{\vh_{\rm mag}}
\newcommand{\HDZ}{\H}
\newcommand{\kmag}{\kappa_{\rm mag}}

\renewcommand{\H}{\mathcal{H}}
\newcommand{\F}{\mathcal{F}}
\newcommand{\dL}{d\mathcal{L}}
\newcommand{\dR}{d\mathcal{R}_g}
\newcommand{\R}{\mathbb{R}}
\newcommand{\T}{\mathbb{T}}
\newcommand{\N}{\mathbb{N}}
\newcommand{\Om}{\Omega}
\newcommand{\om}{\omega}
\newcommand{\la}{\lambda}
\newcommand{\up}{\Upsilon}
\newcommand{\tup}{\tilde{\up}}

\numberwithin{equation}{section}

\begin{document}

\maketitle

\begin{abstract}
The goal of this paper is to describe Zermelo's navigation problem on Riemannian manifolds as a time-optimal control problem and give an efficient method in order to evaluate its control curvature.
We will show that up to change the Riemannian metric on the manifold the control curvature of Zermelo's problem has a simple to handle expression which naturally leads to a generalization of the classical Gauss-Bonnet formula in an inequality. This Gauss-Bonnet inequality enables to generalize for Zermelo's problems the E. Hopf theorem on flatness of Riemannian tori without conjugate points.\\

\end{abstract}


\section{Introduction}

In the present paper we study a special class of time-optimal control problems on two-dimensional manifolds: the Zermelo'-like problems.
By Zermelo'-like problems we mean the class of time-optimal control problems formed by the classical Zermelo's navigation problems on Riemannian manifolds and the corresponding co-problems.

Our first goal in this paper is to describe these two problems and give an explicit expression for their control curvature, the latter being the control analogue of the Gaussian curvature of surfaces. This is the purpose of Section \ref{Z'-like_pb}.  

Zermelo's navigation problem aims to find the minimum time trajectories in a Riemannian manifold $(M,g)$ under the influence of a drift represented by a vector field \(\X\).
The study of Zermelo's navigation problem began in 1931 with the work by E. Zermelo \cite{Zermelo1931} and a while latter by C. Carath{\'e}odory in \cite{cara}.
In a recent paper \cite{BaoRoblesShen2004}, Zermelo's navigation problem has been studied has a special case of Finslerian metrics and has been an efficient tool in order to give a complete classification of strongly convex Randers metrics of constant flag curvature, the latter being the Finslerian analogue of the Riemannian sectional curvature.

The co-Zermelo's navigation problem on a Riemannian surface $(M,g)$ with drift $\up$, where $\up$ is a one-form on $M$, is a time-optimal for which the maximized Hamiltonian function $h$ resulting of Pontryagin Maximum Principle has level one equal to the fiber bundle
$\cup_{q \in M}\{\up_q + \mathcal{S}^{g*}_q\}$, where $\mathcal{S}^{g*}_q$ is the unitary Riemannian cosphere of the metric $g$.
What is surprising with this definition of the problem is that it naturally leads to choose a good system of coordinates in which the control curvature has very nice and simple expression as a function of the drift one-form and the Gaussian curvature of the metric $g$.
Contrary to the Zermelo's navigation problem the curvature of the co-Zermelo problem is much more readable than the one Zermelo problem itself and thus, much more exploitable.

Another surprising property of the co-Zermelo's problem is that its flow is just a time rescaling of the magnetic flow of the pair $(g,d\up)$, the latter being solution of a fixed time variational problem.
In particular, it implies that the curvature of the problem of a charged particle in a magnetic field is just a reparametrization of the curvature of co-Zermelo problem.

We prove constructively that Zermelo's navigation problem on $(M,g)$ with drift vector field $\X$ is feedback equivalent to a co-Zermelo problem but this time with respect to another Riemannian metric. This is the contents of Proposition \ref{EquivZerDZ} and its Corollary \ref{feedback_equiv_Zer-coZer} in \S \ref{duality}.
This proposition is fundamental because it points out that there are two different Riemannian metrics canonically associated to a given Zermelo problem.
In particular, it implies that the two problems have the same curvature and also allows to see a Zermelo's navigation problem as its dual co-Zermelo problem and vice versa.

This is of particular interest because
the presentation of a given Zermelo's navigation problem as its feedback equivalent co-Zermelo problem has the serious advantage to present the curvature of the considered problem in an easier to handle formula.
It is thus, clearly, the way to do.
It also show how the classical Zermelo's navigation problem is linked to magnetic flows.

It is the second goal of this paper to show that there is a natural way to generalize the classical Gauss-Bonnet formula for Riemannian surfaces to an inequality for Zermelo-like problems. More precisely, we will see that, given a Zermelo'-like problem on a surface $M$ there exists a canonically defined positive function $\phi$ such that
$\int_\H \phi \kappa \, \dL \ge \chi(M)$,
where $\chi(M)$ is the Euler characteristic of $M$, $\H$ is the hypersurface $h^{-1}(1)$ and $\dL$ is the Liouville volume on $\H$.
Moreover, the function $\phi$ takes the constant value equal to one if and only if the Zermelo problem is indeed Riemannian, in which case the inequality turns to be the classical Gauss-Bonnet formula.
This Theorem \ref{GB} in \S \ref{SectionGB} of the paper.

Our last goal in this paper is to generalize for Zermelo's problems the E. Hopf's theorem which asserts that two-dimensional Riemannian tori without conjugate points are flat. This will be done in two steps following the Hopf's method.
First we show that if a control system on a compact surface without boundary has no conjugate points then, its total curvature $\int_\H \kappa \, \dL$ must be negative or zero and, in the latter case its curvature must be zero identically.
This is Theorem \ref{E.Hopf_Th}.
The second step is to use the Gauss-Bonnet inequality together with Theorem \ref{E.Hopf_Th} to deduce flatness.
In the Riemannian situation Theorem \ref{E.Hopf_Th} together with the Gauss-Bonnet inequality (which, in this case, reduces to the classical Gauss-Bonnet formula) imply straightforwardly flatness for tori without conjugate points.
Of course Theorem \ref{E.Hopf_Th} applies to Zermelo'-like problems but,
due to the presence of the function $\phi$ in the Gauss-Bonnet inequality, the situation is more delicate and essentially different.
Indeed, a Zermelo'-like problem without conjugate points on a Riemannian torus is not necessarily flat unless its total curvature is zero.
This situation is described in Theorem \ref{E.Hopf_Th_For_Zer} and its Corollaries \ref{E.Hopf_Th_for_DZ_on_T^2} and \ref{Flat_Zer'-like_pbs_on_T^2}.

To conclude our paper we discuss the further generalizations of the presented results to more general situations than the Riemannian one. We will see that even in the special case of Landsberg surfaces surfaces not all results can be transposed.

\section{Curvature of two-dimensional smooth control systems}
In the present paper smooth objects are supposed to be of class $\mathcal{C}^\infty$.
Let us fix some notations.
For a two-dimensional manifold $M$,
$\pi: T^*M \to M$ is the cotangent bundle to $M$.
We denote by $s$ the canonical Liouville one-form on $T^*M$,
$s_\la=\la \circ \pi_*$, $\la \in T^*M$.

If $M$ is endowed with a Riemannian structure $g$, $\lang \cdot , \cdot \rang_g$ and $|\cdot|_g$ denote the Riemannian scalar product and the Riemannian norm respectively. Since the Riemannian structure defines a canonical identification between the tangent and cotangent bundle of $M$, we use the notations of the scalar product and norm indifferently for vectors and covectors, vector fields and one-forms.
We denote by $\mathcal{S}^g$ and $\mathcal{S}^{g*}$
the unitary spherical bundle $\{ v   \in TM   \, : \, |v|_g   = 1 \}$
and the unitary cospherical bundle $\{ \xi \in T^*M \, : \, |\xi|_g = 1 \}$ respectively.

\subsection{Definition}
We briefly recall some facts concerning the curvature of smooth control systems in dimension two. For more details on the subject we refer the reader to one of the following items \cite{AAAbook, ulysseThese, ulysseDEA}.

Consider the following time-optimal smooth control problem 
\begin{eqnarray}
&& \dot{q} = \bs{f}(q,u),\quad q \in M,\quad u \in U, \label{optpb_dyn}\nonumber \\
&& q(0) = q_0, \quad q(t_1) = q_1, \label{optpb_bc} \\
&& t_1 \to \min, \label{optpb_cost} \nonumber
\end{eqnarray}
where $M$ and $U$ are connected smooth manifolds of respective dimension two and one.
For the above time-optimal control problem we denote
by $h = \max_{u \in U}\lang \la,\bs{f}(q,u) \rang$, $\la \in T^*_qM$, $q \in M$, the (normal) Hamiltonian function of PMP (Pontryagin Maximum Principle),
by $\H = h^{-1}(1) \subset T^*M$,
and by $\vh$ the Hamiltonian field associated with the restriction of $h$ to $\H$.
Recall that the maximized Hamiltonian $h$ is a function on the cotangent bundle $T^*M$ one-homogeneous on fibers and non-negative.
Under the regularity assumptions of strong convexity 
%
\begin{equation}\label{Reg_Assumptions}
\bs{f}(q,u) \wedge \parfrac{\bs{f}(q,u)}{u} \neq 0,\quad
\parfrac{\bs{f}(q,u)}{u} \wedge \parfrac{^2\bs{f}(q,u)}{u^2} \neq 0, \quad
q \in M, \quad u \in U,
\end{equation}
the curve $\H_q = \H \cap T^*_qM$ admits, up to sign and translation, a natural parameter providing us with a vector field
$\bs{v}_q$ on $\H_q$ and by consequence with a vertical vector field $\bs{v}$ on $\H$.
Vector field $\bs{v}$ is characterized by the fact that it is, up to sign, the unique vector field on
$\H$ such that
\begin{equation}\label{b}
L^2_{\bs{v}}s|_\H = -s|_\H + b L_{\bs{v}}s|_\H,
\end{equation}
where
$b$ is a smooth
function on the level $\H$.
The function $b$, which is by definition a feedback-invariant, is called the centro-affine curvature.

The vector fields $\vh$ and $\bs{v}$ which are, by definition, feedback-invariant satisfy the following nontrivial commutator relation:
\begin{equation}\label{courbure}
\Big[ \vh , \Big[ \bs{v} , \vh \Big] \Big] = \kappa \bs{v},
\end{equation}
where the coefficient $\kappa$ is defined to be {\it the control curvature} or simply {\it the curvature} of the optimal control problem (\ref{optpb_dyn})-(\ref{optpb_cost}).
The control curvature is by definition  a feedback-invariant of the control system and a function on $\H$ (and not on $M$ as the Gaussian one). Moreover, $\kappa$ is the Gaussian curvature if the control system defines a Riemannian geodesic problem.

\begin{example}\label{Riem_geo_pb}
{\rm
Consider the time-optimal control problem corresponding to the geodesic problem on a two-dimensional Riemannian surface $(M,g)$:
\begin{eqnarray*}
&& \dot{q} = u, \quad q \in M, \quad u \in T_qM, \quad |u|_g = 1, \\
&& q(0) = q_0,\quad q(t_1) = q_1 \\
&& t_1 \to \min, \nonumber
\end{eqnarray*}
In this case, the Hamiltonian function of PMP is given by
\begin{equation*}
h_g(\la) = |\la|_g, \quad \la \in T^*M,
\end{equation*}
and the vectors fields $\vh_g$ and $\bs{v}_g$ on $h_g^{-1}(1)$ by
\begin{equation*}
\vh_g = \cos \theta\,\bs{e}_1+\sin \theta\,\bs{e}_2 + (c_1\cos \theta + c_2\sin \theta)\parfrac{}{\theta},
\quad \bs{v}_g = \parfrac{}{\theta},
\end{equation*}
where $(\bs{e}_1,\bs{e}_2)$ is a local $g$-orthonormal frame whose structural constants
$c_1$, $c_2$ are defined on $M$ by
\begin{equation*}
[\bs{e}_1,\bs{e}_2]=c_1\bs{e}_1+c_2\bs{e}_2,\quad c_1,\ c_2\in C^\infty(M),
\end{equation*}
and $\theta$ is the parameter on the fiber $h_g^{-1}(1) \cap T^*_qM = \mathcal{S}^{g*}_q$ defined by
\begin{equation*}
\lang \la ,\bs{e}_1(q) \rang=\cos\theta,\quad
\lang \la ,\bs{e}_2(q) \rang=\sin\theta.
\end{equation*}
The Gaussian curvature $\kappa_g$ of the surface $(M,g)$ is evaluated as follows:
\begin{equation}\label{riemcurv}
\kappa_g = -c_1^2 -c_2^2 + L_{\bs{e}_1}c_2 - L_{\bs{e}_2}c_1.
\end{equation}
Of course, for the Riemannian problem the curvature depends only on the base point $q\in M$ as one can see from formula (\ref{riemcurv}) but in
general this is not the case: the control curvature depends also on the coordinate in the fiber $\mathcal{H}_q$ 
and thus is a function on the whole three-dimensional manifold $\mathcal{H}$.
}\end{example}
%
\subsection{Reparametrization}\label{rep-curv}
In practice, it is sometimes easier for computations to not to consider the curvature itself but some reparametrization of it.
We will thus see how the curvature $\kappa$ changes under a reparametrization of time.
So let $t = T(\tau)$ be a reparametrization of time. Under this reparametrization the ODE 
$\frac{d\la}{dt} = \vh(\la)$ changes as follows:
\begin{equation*}
\frac{d\la \circ T}{d\tau} = \frac{d\la}{dt}\frac{dT}{d\tau} = \frac{dT}{dt} \vh(\la).
\end{equation*}
Thus, reparametrizing time just means to consider the field $\vh$ in the form
\begin{equation*}\label{rep-h}
\vh = \frac{\bs{\hat{h}}}{\varphi},
\end{equation*}
where $\varphi \in C^\infty(\H)$ is a positive function whose primitive along the trajectories of $\vh$
is the time reparametrization function.
Let $\ve$ be is the vertical field defined by
\(\bs{v} = \sqrt{\varphi}\ve.\)
Vector fields $\he$ and $\ve$ satisfy the following non trivial commutator relation:
\begin{equation}\label{rep-courbure}
\Big[ \he , \Big[ \ve , \he \Big] \Big] = \hat{\kappa}\ve + \xi\he,
\quad \hat{\kappa},\ \xi\in C^\infty(\H).
\end{equation}
Indeed, denote for simplicity $\psi = \sqrt{\varphi}$ then, we have:
\begin{eqnarray*}
\Big[ \he , \Big[ \ve , \he \Big] \Big]
&=&  \Big[ \psi^2\vh , \Big[ \psi^{-1}\bs{v} , \psi^2\vh \Big] \Big]
 = \Big[ \psi^2\vh ,
         \psi\Big[ \bs{v} , \vh \Big] + \psi^{-1}L_{\bs{v}}(\psi^2)\vh - \psi^2L_{\vh}(\psi^{-1})\bs{v} \Big] \\
&=&  \Big[ \psi^2\vh , \psi \Big[\bs{v} , \vh \Big] \Big] 
   + \Big[ \psi^2\vh , L_{\vh}\psi\bs{v}\Big] + \xi_1\vh \\
&=&  \psi^3 \Big[ \vh , \Big[ \bs{v} , \vh \Big] \Big]
   + \psi^2 L_{\vh}\psi \Big[ \bs{v} , \vh \Big]
   + \psi^2 L_{\vh}\psi \Big[ \vh , \bs{v} \Big]
   + \psi^2 L_{\vh}^2 \psi \bs{v} + \xi_2\vh \\
&=&  (\psi^4\kappa+\psi^3L_{\vh}^2\psi)\ve+\xi\he \\
&=&  \hat{\kappa} \ve + \xi\he,
\end{eqnarray*}
where $\xi_1$, $\xi_2$, $\xi \in C^\infty(\H)$.
From the previous computation one infers that the curvature and its reparametrization $\hat{\kappa}$ satisfy the following relation:
\begin{equation}\label{K=Krho-S}
\kappa = \frac{\hat{\kappa} - \mathcal{S}(\varphi)}{\varphi^2}, \quad
\mathcal{S}(\varphi) = \varphi L_{\vh}\left( \frac{L_{\vh}\varphi}{2} \right)
                      -\left( \frac{L_{\vh}\varphi}{2} \right)^2.
\end{equation}
%

We call the function $\hat{\kappa}$ defined by the relation (\ref{rep-courbure}) the
{\it $\varphi$-reparametrization} of the curvature $\kappa$.

\section{Zermelo'-like problems}\label{Z'-like_pb}
\subsection{Zermelo's navigation problem}
In his article \cite{Zermelo1931} of 1931 Ernst Zermelo formulates the following problem:\\

\noindent{\textit{``In an unbounded plane where the wind distribution is given by a vector field as
a function of position and time, a ship moves with constant velocity relative to the
surrounding air mass. How must the ship be steered in order to come from a starting point
to a given goal in the shortest time?"}}\\

For our purpose we assume that we are working on a Riemannian surface in the presence of a stationary wind distribution that we call drift.
Zermelo's navigation problem thus consists of finding the quickest path (in time) of a point on a Riemannian surface $(M,g)$ in the presence of a stationary drift modeled by an autonomous vector field $\bs{X}\in\Vec{M}$.
This time-optimal control problem is read
\begin{eqnarray}
&& \dot{q}=\bs{X}(q)+u,\quad q\in M,\quad u\in T_qM, \quad |u|_g = 1, \label{Zer_dyn} \\
&& q(0)=q_0,\quad q(t_1)=q_1                                          \label{Zer_bc} \\
&& t_1\to \min,                                                       \label{Zer_cost}
\end{eqnarray}
and we call it {\it Zermelo problem of the pair $(g,\X)$}.
The Hamiltonian function of PMP is
\begin{equation}\label{hamZer}
\hZ(\la)=
\max_{|u|_g \le 1} \left( \lang \la , \X \rang + \lang \la , u \rang \right)
=\lang \la,\X(q)\rang + |\la|_g,
\end{equation}
and the Hamiltonian vector field on $\H = \hZ^{-1}(1)$ has the form
\begin{equation}\label{ham_vec_zermelo}
\vhZ
= \X+\vh_g
+ \Big( \lang u_{\rm max} , [\bs{e}_1,\bs{e}_2]\rang_g \lang u_{\rm max} , \X \rang_g
+ L_{[u_{\rm max},\bs{v}_g]}\lang u_{\rm max} , \X \rang_g \Big) \bs{v}_g,
\end{equation}
where the function $u_{\max} = u_{\max}(\la)$ is the restriction to $\H$ of the maximized control obtained in the maximization
(\ref{hamZer}).
Relation (\ref{ham_vec_zermelo}) leads naturally to an expression of the curvature of Zermelo's navigation problem (\ref{Zer_dyn})-(\ref{Zer_cost}) as a function of the drift $\X$ and the Gaussian curvature of the surface $(M,g)$. 
We do not give here a precise formula for this expression of the curvature since it leads to a formula which is rather complicated and hardly exploitable except for very simple cases.
We refer the reader to
\cite{ulysseDEA} for a detailed description and coordinate expression of the curvature of this problem.

\subsection{Co-Zermelo's navigation problem}\label{SectionDZ}
Roughly speaking, whereas Zermelo's navigation problem was defined by its dynamics,
i.e., as a subbundle of the tangent bundle over the state space $M$
(in this case $\cup_{q \in M}\{ \X(q) + \mathcal{S}^g_q \} \subset TM$), 
co-Zermelo's navigation problem will be defined as a subbundle of $T^*M$.
Precisely,
\begin{definition}
We call co-Zermelo problem of the pair $(g,\up)$ the minimum time problem for which the Hamiltonian function of PMP has level one equals to
$\cup_{q \in M} \{ \up_q + \mathcal{S}^{g*}_q \} \subset T^*M$
where $\up$ is a one-form on $M$ such that $|\up|_g < 1$.
\end{definition}
Let $\hDZ$ be the maximized Hamiltonian function of PMP associated to the co-Zermelo problem of the pair $(g,\up)$ which, let us recall it, is
one-homogeneous on fibers and non-negative.
Denote by $\HDZ$ the hypersurface $\hDZ^{-1}(1)$.
By definition of the co-Zermelo problem the hypersurface $\HDZ$ is characterized by
\begin{equation}\label{DZdef}
\lang \lambda - \up_{\pi(\la)} , \lambda - \up_{\pi(\la)}\rang_g = 1, \quad \forall \lambda \in \HDZ.
\end{equation}
Suppose now that $\la \in T^*M$ is a non zero covector such that $\hDZ(\la)\neq 0$. Then, using the homogeneity of $\hDZ$ we get
\begin{equation*}
\frac{\la}{\hDZ(\la)}\in \HDZ.
\end{equation*}
Consequently, the covector $\la/\hDZ(\la)$ has to satisfy equation (\ref{DZdef}). Plugging this covector in equation (\ref{DZdef}) leads to
\begin{equation}\label{hameq}
\lang \la - \hDZ(\la)\up_{\pi(\la)} , \la - \hDZ(\la)\up_{\pi(\la)}\rang_g^2 = \hDZ(\la)^2, \quad \la \in T^*M,
\end{equation}
which gives an implicit definition for the Hamiltonian function $\hDZ$. Solving equation (\ref{hameq}) for $\hDZ(\la)$ gives
\begin{equation}\label{Ham_func_coZer}
\hDZ(\la)
= \frac{ - \lang \la , \up_{\pi(\la)} \rang_g 
         + \sqrt{\lang \la , \up_{\pi(\la)} \rang_g^2 + \left( 1 - |\up_{\pi(\la)}|_g^2 \right) |\la|_g^2}}
{1 - |\up_{\pi(\la)}|_g^2},
\end{equation}
where we have excluded the non-positive solution.
We now derive the equation of the Hamiltonian field associated to $\hDZ$ on the level surface $\HDZ$.
If $(p,q)$ is a canonical system of local coordinates on $T^*M$ and $\theta$ is the coordinate on fibers
$\HDZ \cap T^*_qM$, the vector field $\vhDZ$ is given by
\begin{equation*}
\vhDZ = \pi_*\vhDZ + c\parfrac{}{\theta}
= \left. \parfrac{h}{p_1} \right|_{\HDZ} \parfrac{}{q_1} + \left. \parfrac{h}{p_2} \right|_{\HDZ} \parfrac{}{q_2}
 + c\parfrac{}{\theta}.
\end{equation*}
We define a parameter $\theta$ on fibers $\HDZ \cap T^*_{\pi(\la)}M$ in the following manner.
Let $(\bs{e}_1,\bs{e}_2)$ be a local $g$-orthonormal frame on $M$.
Notice that equation (\ref{DZdef}) can be locally rewritten
\begin{equation*}
\lang \lambda - \up_{\pi(\la)},\bs{e}_1(\pi(\la)) \rang^2 + \lang \lambda - \up_{\pi(\la)},\bs{e}_2(\pi(\la)) \rang^2=1,
\quad \forall \lambda\in\HDZ.
\end{equation*}
Hence, the fiber $\H_q$ can be naturally parametrized by an angle $\theta$:
\begin{equation*}\label{(theta,q)DualZer}
\lang \la - \up_{\pi(\la)},\bs{e}_1(\pi(\la)) \rang = \cos\theta,\quad
\lang \la - \up_{\pi(\la)},\bs{e}_2(\pi(\la)) \rang = \sin\theta.
\end{equation*}
In order to get the equations of the Hamiltonian vector field
$\vhDZ$, we write equation (\ref{hameq}) in coordinates $(p,q)$
\begin{equation*}
\lang p - \hDZ(p,q)\up_q , \bs{e}_1(q) \rang^2 + \lang p - \hDZ(p,q)\up_q , \bs{e}_2(q) \rang^2
= \hDZ^2(p,q),
\end{equation*}
and we differentiate it with respect to the $p_i$'s. We get
\begin{equation*}
\sum_{k=1}^{2} \lang p - \hDZ\up,\bs{e}_k \rang \lang \parfrac{p}{p_i} - \parfrac{\hDZ}{p_i}\up,\bs{e}_k \rang
= \hDZ\parfrac{\hDZ}{p_2},
\quad i = 1,2.
\end{equation*}
Consequently, on the surface $\HDZ$
\begin{eqnarray*}
\left. \parfrac{\hDZ}{p_i} \right|_{\HDZ}
&=& \frac{\lang p - \up , \bs{e}_1 \rang e_1^i + \lang p - \up,\bs{e}_2 \rang e_2^i}
{
1
+ \lang p-\up , \bs{e}_1 \rang \lang \up , \bs{e}_1 \rang
+ \lang p-\up , \bs{e}_2 \rang \lang \up , \bs{e}_2 \rang
} \\
&=& \frac{\cos\theta\,e_1^i+\sin\theta\,e_2^i}
{1 + \cos\theta \lang \up , \bs{e}_1 \rang + \sin\theta \lang \up , \bs{e}_2 \rang},
\quad i = 1,2.
\end{eqnarray*}
Thus the horizontal part of the field $\vhDZ$ on $\HDZ$ is
\begin{eqnarray*}
\pi_*\vhDZ
&=&\frac{1}{\phicoZ}
\left(
(\cos\theta\,e_1^1+\sin\theta\,e_2^1)\parfrac{}{q_1}+
(\cos\theta\,e_1^2+\sin\theta\,e_2^2)\parfrac{}{q_2}
\right) \\
&=& \frac{1}{\phicoZ}(\cos\theta\,\bs{e}_1+\sin\theta\,\bs{e}_2)
\end{eqnarray*}
where
\begin{equation*}
\phicoZ(\theta,q)
= 1 + \cos\theta \lang \up_q , \bs{e}_1(q) \rang + \sin\theta \lang \up_q , \bs{e}_2(q) \rang.
\end{equation*}
Because $\vh$ is the Hamiltonian field in restriction to $\H$,
we have $ds|_\H(\vh,\cdot) = 0$
from which we can deduce the $\parfrac{}{\theta}$ of $\vh$.
Let $(\bs{e}_1^*,\bs{e}_2^*)$ be the coframe dual to $(\bs{e}_1,\bs{e}_2)$ and denote $s|_\H = \om$. 
In coordinates $\la = (\theta,q)$ on $\H$ the Liouville one-form $\om$ takes the form
%
\begin{eqnarray}
\om 
&=& \lang \la,\bs{e}_1 \rang \bs{e}_1^* + \lang \la,\bs{e}_2 \rang \bs{e}_2^*  \nonumber\\
&=& \left( \lang \la - \up,\bs{e}_1 \rang + \lang \up,\bs{e}_1 \rang \right) \bs{e}_1^*
   +\left( \lang \la - \up,\bs{e}_2 \rang + \lang \up,\bs{e}_2 \rang \right) \bs{e}_2^*  \nonumber\\
&=& \cos\theta\bs{e}_1^* + \sin\theta\bs{e}_2^* + \up, \label{LiouOnHDualZer}
\end{eqnarray}
so that its exterior derivative is
\begin{equation*}
d\om
= -\sin\theta d\theta\wedge\bs{e}_1^* + \cos\theta d\theta\wedge\bs{e}_2^*
  +\cos\theta d\bs{e}_1^* + \sin\theta d\bs{e}_2^* + d\up.
\end{equation*}
Using Cartan's formula for one forms
$d\xi(\bs{X},\bs{Y}) = L_{\X}\lang \xi,\bs{Y} \rang - L_{\bs{Y}}\lang \xi,\bs{X} \rang - \langle \xi,[\X,\bs{Y}] \rangle$,
one easily see that
\begin{equation*}
d\bs{e}_1^* = -c_1 dV_g, \quad d\bs{e}_2^* = -c_2 dV_g,
\end{equation*}
where, as in Section \ref{Riem_geo_pb}, $c_1$, $c_2$, are the structural constants of the frame $(\bs{e}_1,\bs{e}_2)$ and,
$dV_g = \bs{e}_1^* \wedge \bs{e}_2^*$ denotes the Riemannian volume element on M.
Let $\Om \in \mathcal{C}^\infty(M)$ be the function defined by $d\up=-\Om\, dV_g$
and denote $c_g = c_1 \cos\theta + c_2 \sin\theta$.
Summing up, we have
\begin{equation}\label{dw_coZer}
d\om
= -\sin\theta d\theta\wedge\bs{e}_1^* + \cos\theta d\theta\wedge\bs{e}_2^*
  -(c_g + \Om)dV_g,
\end{equation}
from which we get
\begin{equation*}
0 = d\om(\vh, \cdot)
  = - c \sin\theta\bs{e}_1^* + c \cos\theta\bs{e}_2^* 
  + \frac{c_g+\Om}{\phicoZ}\sin\theta\bs{e}_1^* - \frac{c_g+\Om}{\phicoZ}\cos\theta\bs{e}_2^*.
\end{equation*}
Hence,
\begin{equation*}
c = \frac{c_g+\Om}{\phicoZ}.
\end{equation*}
Summing up, the Hamiltonian of the co-Zermelo problem reads
\begin{equation*}
\vhDZ(\theta,q)
= \frac{1}{\phicoZ(\theta,q)} \left( \cos\theta\,\bs{e}_1(q) + \sin\theta\,\bs{e}_2(q) 
                                    +(c_g(q) + \Om(q))\parfrac{}{\theta}
                              \right)
\end{equation*}
or, equivalently
\begin{equation}\label{hriemgvriem}
\vhDZ
= \frac{1}{\phicoZ}\left( \diffeo{}_*\vh_g+\Om \diffeo{}_*\bs{v}_g \right), \quad
\phicoZ(\la) = 1 + \lang \la - \up_{\pi(\la)} , \up_{\pi(\la)} \rang_g, \quad \la \in \H,
\end{equation}
where $\vh_g$ and $\bs{v}_g$ are defined as in Section \ref{Riem_geo_pb}
and $\diffeo{}$ is the diffeomorphism
\begin{equation}\label{diffeo_F}
\begin{array}{rcl}
\diffeo{} : \mathcal{S}^{g*} & \to & \H \\
\la & \mapsto & \la +\up_{\pi(\la)}.
\end{array}
\end{equation}
Notice that $(\diffeo{})^{-1} = \diffeo{-}$.

\begin{remark}
{\rm
To conclude this section let us give a (coordinate free) formulation for the co-Zermelo problem of the pair $(g,\up)$ as a time-optimal control problem.
According to (\ref{hriemgvriem}), this time-optimal control problem reads
\begin{eqnarray}
&& \dot{q}={\displaystyle\frac{u}{1 + \lang \up_q , u \rang}},\quad q\in M,\quad u\in T_qM,\quad |u|_g=1,\label{coZer_dyn}\\
&& q(0)=q_0,\quad q(t_1)=q_1 \nonumber\label{coZer_bc}\\
&& t_1\to \min,\nonumber\label{coZer_cost}
\end{eqnarray}
and the reader can check that the result of the maximality condition of PMP,
$ \max_{|u|_g = 1} \langle \la , \dot{q} \rangle$,
is the Hamiltonian function given by relation (\ref{Ham_func_coZer}).
}
\end{remark}

\subsection{Curvature of the co-Zermelo problem}
In order to get the expression of the curvature of the co-Zermelo problem, we first of all need to find the expression of the vertical field that satisfies relation (\ref{b}).

According to (\ref{LiouOnHDualZer}) and (\ref{dw_coZer}),
\begin{equation}\label{om^d_qom=phidV}
\om \wedge d\om = \om\wedge \parfrac{\om}{\theta} = \phicoZ\,\bs{e}_1^*\wedge\bs{e}_2^* = \phicoZ\,dV_g \neq 0,
\end{equation}
which shows that
$(\om , \parfrac{\om}{\theta})$ forms a frame of horizontal one-forms on $\HDZ$.
The decomposition of the second derivative $\parfrac{^2\om}{\theta^2}$ in this frame reads
\begin{equation*}
\parfrac{^2\om}{\theta^2}
= - \frac{1}{\phicoZ}\om + \frac{\parfrac{\phicoZ}{\theta}}{\phicoZ}\parfrac{\om}{\theta},
\end{equation*}
from which we deduce that the vertical vector field $\bs{v}$ that satisfies (\ref{b}) has the coordinate expression
\begin{equation}\label{VertVDualZer}
\bs{v} = \sqrt{\phicoZ}\parfrac{}{\theta}.
\end{equation}
We now compute the curvature of the co-Zermelo problem according to relation (\ref{courbure}). We find that
\begin{proposition}
The curvature of the co-Zermelo problem of the pair $(g,\up)$ is
\begin{equation}\label{CourbureDualZer}
\kappa_\coZ^{\sr (g,\up)}
=  \frac{1}{( \phicoZ )^2}\left( \kappa_g + \Om^2 + L_{\diffeo{}_*\left[\vh_g,\bs{v}_g\right]}\Om \right) \circ \pi
  -\frac{\mathcal{S}(\phicoZ)}{(\phicoZ)^2}.
\end{equation}
\end{proposition}
\begin{proof}
According to (\ref{hriemgvriem}) and (\ref{VertVDualZer}),
\begin{displaymath}
\begin{array}{rclcrcl}
\vhDZ    & = & \displaystyle{ \frac{\he}{\phicoZ} }, & \quad & \he & = & \diffeo{}_*\big( \vh_g + \Om\bs{v}_g \big), \\
\bs{v}   & = & \displaystyle{ \sqrt{\phicoZ} }\ve  , & \quad & \ve & = & \diffeo{}_*\bs{v}_g,
\end{array}
\end{displaymath}
which implies that it is enough for this problem to compute the $\phicoZ$-reparametrized curvature (defined in Section \ref{rep-curv}). We have
\begin{eqnarray*}
\diffeo{-}_*\Big[\he,\left[\ve,\he\right]\Big]
&=&  \Big[ \vh_g + \Om\bs{v}_g , \Big[ \bs{v}_g , \vh_g + \Om\bs{v}_g \Big] \Big] \\
&=&  \Big[ \vh_g , \Big[ \bs{v}_g , \vh_g \Big] \Big]
   + \Om \Big[ \bs{v}_g , \Big[ \bs{v}_g , \vh_g \Big] \Big]
   + L_{\left[ \vh_g , \bs{v}_g \right]} \Om\bs{v}_g \\
&=&  \kappa_g\bs{v}_g
   - \Om\vh_g
   + L_{\left[ \vh_g , \bs{v}_g \right]} \Om\bs{v}_g \\
&=& (\kappa_g + \Om^2 + L_{\left[ \vh_g , \bs{v}_g \right]} \Om)\diffeo{-}_*\ve - \Om \diffeo{-}_*\he.
\end{eqnarray*}
According to (\ref{K=Krho-S}) the result follows.
\end{proof}\\

We refer the reader to \cite{piotr} for a detailed presentation of the co-Zermelo problem with linear drift on the Euclidean plane $\R^2$. In particular, using the reparametrized curvature, the author studied in great details the occurrence of conjugate points.

\subsection{Duality between Zermelo and co-Zermelo problems}\label{duality}
In this section we prove a proposition which asserts the feedback equivalence between the Zermelo and the co-Zermelo navigation problems.
Although this proposition is simple indeed, it will have a fundamental role in the sequel due to fact that the curvature is much simpler to handle for the co-Zermelo problem than for the Zermelo navigation problem itself.

Let $(M,g)$ be a Riemannian manifold and fix an $g$-orthonormal frame $(\bs{e}_1,\bs{e}_2)$.

If $\X \in \Vec{M}$, we define the local orthonormal frame for $g$ associated to the vector field $\X$ with respect to the frame $(\bs{e}_1,\bs{e}_2)$ by
\begin{equation*}
\begin{array}{rcrcrl}
\bs{e}_1^{\sr\X} &=&  \cos\theta^{\sr\X}\bs{e}_1 &+& \sin\theta^{\sr\X}\bs{e}_2 &\phantom{\Big(} \\
\bs{e}_2^{\sr\X} &=& -\sin\theta^{\sr\X}\bs{e}_1 &+& \cos\theta^{\sr\X}\bs{e}_2 &,
\end{array}
\end{equation*}
where $q \mapsto \theta^{\scriptscriptstyle \X}(q)$ is the angle defined by
\begin{equation}\label{def_angle_theta_X}
\begin{cases}
\theta^{\sr \X}(q) =0 & \text{if $\X(q) = 0$}, \\
{\displaystyle
\cos\theta^{\sr \X}(q)=\frac{\lang \bs{X}(q),\bs{e}_1(q)\rang_g}{|\bs{X}(q)|_g},
\quad
\sin\theta^{\sr \X}(q)=\frac{\lang \bs{X}(q),\bs{e}_2(q)\rang_g}{|\bs{X}(q)|_g}}
& \text{if $\X(q) \neq 0$}.
\end{cases}
\end{equation}
In the same way if $\up \in \Lambda^1(M)$ we define the $g$-orthonormal frame associated to the one-form $\up$ with respect to the frame $(\bs{e}_1,\bs{e}_2)$ by
\begin{equation*}
\begin{array}{rcrcrl}
\bs{e}_1^{\sr\up} &=&  \cos\theta^{\sr\up}\bs{e}_1 &+& \sin\theta^{\sr\up}\bs{e}_2 &\phantom{\Big(} \\
\bs{e}_2^{\sr\up} &=& -\sin\theta^{\sr\up}\bs{e}_1 &+& \cos\theta^{\sr\up}\bs{e}_2 &,
\end{array}
\end{equation*}
where $q \mapsto \theta^{\scriptscriptstyle \up}(q)$ is the angle defined by
\begin{equation*}
\begin{cases}
\theta^{\sr \up}(q) =0 & \text{if $\up_q = 0$}, \\
{\displaystyle
\cos\theta^{\sr \up}(q)=\frac{\lang \up_q,\bs{e}_1(q)\rang}{|\up_q|_g},
\quad
\sin\theta^{\sr \up}(q)=\frac{\lang \up_q,\bs{e}_2(q)\rang}{|\up_q|_g}}
& \text{if $\up_q \neq 0$}.
\end{cases}
\end{equation*}
Notice that in this frames
\begin{equation*}
\X  = \lang \X  , \bs{e}_1^{\sr\X} \rang_g \bs{e}_1^{\sr\X} = |\X|_g \bs{e}_1^{\sr\X},\quad
\up = \lang \up , \bs{e}_1^{\sr\up} \rang \bs{e}_1^{{\sr\up} *} = |\up|_g \bs{e}_1^{{\sr\up} *}.
\end{equation*}
Suppose for now that the Riemannian norm of the drift in our Zermelo navigation is strictly smaller than one.
\begin{proposition}\label{EquivZerDZ}
Let $(M,g)$ be a Riemannian surface. Let $\X$ be a vector field on $M$ (respectively, $\up$ a one-form on $M$).
There exists on $M$ a new Riemannian metric $\tg = \tg(g,\X)$ (respectively $\tg = \tg(g,\up)$)
and a one-form $\tilde{\up}$ (respectively, a vector field $\bs{\tilde{X}}$) such that the Zermelo problem of the pair
$(g, \X)$ (respectively, the co-Zermelo problem of the pair $(g, \up)$)
and the co-Zermelo problem
of the pair $(\tg, \tilde{\up})$ (respectively, the Zermelo problem of the pair $(\tg, \bs{\tilde{X}})$)
have the same Hamiltonians.
\end{proposition}
\begin{proof}
Consider Zermelo's navigation problem (\ref{Zer_dyn})-(\ref{Zer_cost}) and
let $(\bs{e}_1,\bs{e}_2)$ be an orthonormal frame for the metric $g$.
Define some polar coordinates $(\rho,\theta)$ on the fiber $T^*_qM$ by
\begin{equation*}
\rho = |\la|_g, \quad
\lang \la , \bs{e}_1 \rang = \rho \cos\theta, \quad
\lang \la , \bs{e}_2 \rang = \rho \sin\theta,
\end{equation*}
so that the Hamiltonian (\ref{hamZer}) takes the form 
\begin{equation*}
h(\rho,\theta,q) = \rho\left(|\bs{X}(q)|_g\cos(\theta-\theta^{\sr \X}(q))+1\right),
\end{equation*}
where $\theta^{\sr\X}(q)$ is the angle defined by (\ref{def_angle_theta_X}).
Thus, the curve $\H_q = h^{-1}(1) \cap T^*_qM$ has the polar equation
\begin{equation}
\rho(\theta)=\frac{1}{|\bs{X}(q)|_g\cos(\theta-\theta^{\sr \X}(q))+1}.
\end{equation}
Since $|\bs{X}|_g < 1$, the curve $\H_q$ is an ellipse centered at a focus. Moreover, this ellipse has for $g$ a focal distance
$ c = (\rho(\pi + \theta^{\sr \X}) - \rho(\theta^{\sr \X}))/2 = |\X|_g(1-|\X|_g^2)^{-1} $,
a semimajor distance
$a = \rho(\theta^{\sr \X}) + \rho(\pi + \theta^{\sr \X}) = (1-|\X|_g^2)^{-1}$,
and a semiminor distance
$ b = \sqrt{a^2 - c^2} = (1-|\X|_g^2)^{-1/2} $.

In order to transform Zermelo navigation problem in a co-Zermelo problem, we consider the curve $\H_q$ as the drifted Riemannian cosphere at point $q$ for a new Riemannian structure $\tg$ on the manifold. In other words, we ask the one-forms
\begin{equation*}
\bs{\tilde{e}}^*_1 = \frac{1}{1-|\bs{X}|_g^2} \bs{e}_1^{{\sr\X}*}, \quad
\bs{\tilde{e}}^*_2 = \frac{1}{\sqrt{1-|\bs{X}|_g^2}} \bs{e}_2^{{\sr\X}*} 
\end{equation*}
to form an orthonormal coframe for the new Riemannian structure $\tg$ on the manifold and the one-form
\begin{equation*}
\tup = - c\,  \bs{e}_1^{{\sr\X}*} = - \frac{|\X|_g}{1-|\bs{X}|_g^2} \bs{e}_1^{{\sr\X}*}
\end{equation*}
to be the drift one-form of the co-Zermelo problem on $(M,\tg)$.
The corresponding (new) orthonormal frame $(\bs{\tilde{e}}_1,\bs{\tilde{e}}_2)$ is characterized by
\begin{equation*}
\lang (\bs{\tilde{e}}^*_1,\bs{\tilde{e}}^*_2),(\bs{\tilde{e}}_1,\bs{\tilde{e}}_2)\rang={\rm Id},
\end{equation*}
which leads to
\begin{equation*}
\bs{\tilde{e}}_1 = \left( 1-|\bs{X}|_g^2 \right) \bs{e}_1^{\sr\X}, \quad
\bs{\tilde{e}}_2 = \sqrt{1-|\bs{X}|_g^2} \, \bs{e}_2^{\sr\X}.
\end{equation*}
Notice that we have
$( \bs{\tilde{e}}_1 , \bs{\tilde{e}}_2 ) = ( \bs{\tilde{e}}_1^{\sr-\tup} , \bs{\tilde{e}}_2^{\sr-\tup} )$
which shows in particular that
$|\X|_g = |\tup|_{\tg}$.
%

The situation discribed above is illustrated by the picture below.
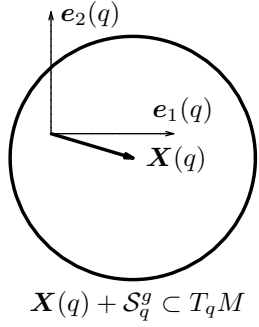
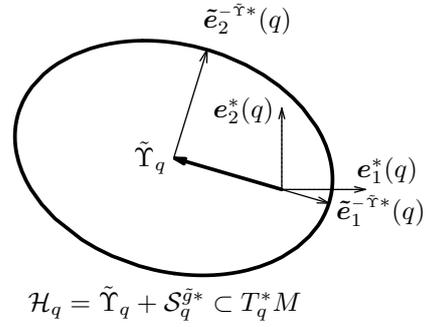
\begin{figure}[!h]
\begin{center}
\input{Fig1EquivZerDZ.pstex_t}
\end{center}
\end{figure}

In order to complete the proof it remains to check that the Hamiltonian function $h_\Z^{\sr (g,\X)}$ of the Zermelo problem of the pair $(g,\X)$ and the Hamiltonian function $h_\coZ^{\sr (\tg,\tup)}$ of the co-Zermelo problem of the pair $(\tg,\tup)$ are the same.
For simplicity we denote $\tilde{c} = |\tup|_{\tg} = |\X|_g$.
We have
\begin{eqnarray*}
h_\Z^{\sr (g,\X)}(\la)
&=&  \lang \la , \X \rang + |\la|_g
 =   \lang \la , \tilde{c} \bs{e}_1^{\sr\X} \rang
   + \sqrt{\lang \la , \bs{e}_1^{\sr\X} \rang^2 + \lang \la , \bs{e}_2^{\sr\X} \rang^2} \\
&=&  \lang \la , \tilde{c} \frac{\bs{\tilde{e}}_1}{1-\tilde{c}^2} \rang
   + \sqrt{\lang \la , \frac{\bs{\tilde{e}}_1}{1-\tilde{c}^2} \rang^2 + \bigg\langle \la , \frac{\bs{\tilde{e}}_2}{\sqrt{1-\tilde{c}^2}} \bigg\rangle^2} \\
&=&
\frac{
\lang \la , \tilde{c} \bs{\tilde{e}}_1 \rang + \sqrt{\lang \la , \bs{\tilde{e}}_1 \rang^2 + (1-\tilde{c}^2)\lang \la , \bs{\tilde{e}}_2 \rang^2}
}{1-\tilde{c}^2}\\
&=&
\frac{
\lang \la , \tilde{c} \bs{\tilde{e}}_1 \rang
+ \sqrt{\lang
\la , \bs{\tilde{e}}_1 \rang^2 + \lang \la , \bs{\tilde{e}}_2 \rang^2 - \tilde{c}^2\lang \la , \bs{\tilde{e}}_2 \rang^2 - \tilde{c}^2\lang \la , \bs{\tilde{e}}_1 \rang^2
+ \tilde{c}^2\lang \la , \bs{\tilde{e}}_1
\rang^2}
}{1-\tilde{c}^2} \\
&=&
\frac{
-\lang \la , -\tilde{c} \bs{\tilde{e}}_1 \rang
+ \sqrt{
\big( \lang \la , \bs{\tilde{e}}_1 \rang^2 + \lang \la , \bs{\tilde{e}}_2 \rang^2 \big)( 1 - \tilde{c}^2 ) 
+ \big( -\tilde{c} \lang \la , \bs{\tilde{e}}_1 \rang \big)^2
}
}{1-\tilde{c}^2} \\
&=&
\frac{
- \langle \la , \langle \tup , \bs{\tilde{e}}_1 \rangle \bs{\tilde{e}}_1^* \rangle_{\tg}
+ \sqrt{|\la|_{\tg}\left( 1- \tilde{c}^2 \right) + \big( \langle \tup , \bs{\tilde{e}}_1 \rangle \lang \la , \bs{\tilde{e}}_1 \rang \big)^2}
}{1-\tilde{c}^2} \\
&=& \frac{
- \langle \la , \tup \rangle_{\tg}
+ \sqrt{ \big( 1 - |\tup|_{\tg}^2 \big) |\la|_{\tg} + \langle \la , \tup \rangle_{\tg}^2 }
}{1 - |\tup|_{\tg}^2} \\
&=& h_\coZ^{\sr (g,\tup)}(\la).
\end{eqnarray*}
In order to prove the converse, one has just to permute the roles of vector fields and one forms in the previous considerations.
\end{proof}\\

Zermelo's navigation problem and co-Zermelo's navigation problem which have the same Hamiltonian are said to be {\it dual problems}.
The above proposition implies in particular that the two dual problems have the same curvature.
This proposition can be reformulated as follows.

\begin{coro}\label{feedback_equiv_Zer-coZer}
Two dual Zermelo's problems are feedback equivalent.
\end{coro}
\begin{proof}
Notations are these of the proof of the previous proposition.
A similar computation computation as the one made in the previous proof shows that the two dual Zermelo's problems have the same sets of admissible velocities, i.e., that for every $q \in M$,
$\{ \X(q)+u \,:\, u \in \mathcal{S}^g_q \} = \{ \tilde{u}(1 - \langle \tup_q,\tilde{u} \rangle)^{-1} \,:\, \tilde{u} \in \mathcal{S}^{\tg}_q \}$ (refer to equations (\ref{Zer_dyn}) and (\ref{coZer_dyn}) for the dynamics of Zermelo's problems).
Thus, the feedback transformation $u \mapsto \tilde{u}(1 - \langle \tup_q,\tilde{u} \rangle)^{-1}- \X(q)$ has obviously the required properties.
\end{proof}
%

\subsection{Classical particle in a magnetic field on a Riemannian surface}
The motion of a charged particle of unit mass under the presence of a magnetic field is modeled by what is called the magnetic flow.
We will see here how the problem of a charged particle in a magnetic field is linked to the dual to Zermelo problem.
Magnetic flows were first considered by Arnold in \cite{Arnold1961} and by Anosov and Sinai in \cite{AnosovSinai1967}
but, it is  Sternberg in \cite{Sternberg1977} gave the first formulation of this problem using symplectic geometry.

Let $(M,g)$ be a two-dimensional Riemannian manifold and $B \in \Lambda^2(M)$ a closed two-form
thought as a magnetic field in which we have absorbed the electric charge of the particle as a parameter.

The magnetic flow of the pair $(g,B)$ is the flow of the Hamiltonian
\begin{equation*}
h_g^2(\la) = \lang \la , \la \rang_g,
\end{equation*}
with respect to the symplectic form $\sigma_B = ds + \pi^*B$ (see \cite{Sternberg1977}).
In the case where $B$ derives from a magnetic potentiel, i.e., when $B = d\up$, $\up \in \Lambda^1(M)$, the magnetic flow is also Hamiltonian with respect to the canonical symplectic form $ds$ but this time with the Hamiltonian function
\begin{equation*}
h_{\rm mag}(\la)  = \frac{1}{2} \lang\la-\up_{\pi(\la)},\la-\up_{\pi(\la)}\rang_g= \frac{1}{2} h_g^2(\la-\up_{\pi(\la)}).
\end{equation*}
A straightforward computation shows that the Hamiltonian vector field $\vh_{\rm mag}$ associated to $h_{\rm mag}$ in restriction to $h_{\rm mag}^{-1}(1)$ is given by
\begin{equation*}
\vh_{\rm mag} = \diffeo{}_*(\vh_g + \Omega\bs{v}_g),
\end{equation*}
where $\Om \in C^\infty(M)$ is defined in same way as the function $\Om$ of the co-Zermelo problem. This shows that the equations of motion of a particle in a magnetic field are in fact the equations of motion of the reparametrized co-Zermelo problem.
For this reason we define the curvature \(\kmag^{\sr (g,d\up)}\) of the magnetic flow to be
the $\phicoZ$-reparametrized curvature of the co-Zermelo problem, i.e.,
\begin{equation}\label{kmag_def}
\kmag^{\sr (g,d\up)} = \kappa_g + \Om^2 + L_{\left[ \vh_g , \bs{v}_g \right]} \Om,
\end{equation}
so that,
\begin{equation}\label{kmag=phi-rep_coZ}
\kappa_{\coZ}^{\sr (g,\up)} = (\phicoZ)^{-2}\left( \kmag^{\sr (g,d\up)} - \mathcal{S}(\phicoZ) \right).
\end{equation}
\begin{remark}{\rm
There is a theory on the reduction of the curvature of Hamiltonian flows by first integrals, see \cite{AAAChtchZel2005}.
The reader can check that, what we have defined to be the curvature of the magnetic flow corresponds to the reduced curvature of the Hamiltonian \(h_{\rm mag}\) on the level \(h_{\rm mag}^{-1}(1)\).
}\end{remark}

\section{A Gauss-Bonnet inequality for Zermelo's problems}\label{SectionGB}
This section is dedicated to some global ``Gauss-Bonnet properties" of Zermelo's problems;
key ingredients to prove Hopf's theorem for Zermelo problems (purpose of the next section).

%
On the three-dimensional surface $\H$ there exists a canonical volume element, called Liouville volume element, defined by
$\dL = - s|_\H \wedge ds|_\H$. Since the Liouville one-form $s|_\H$ is invariant by $\vh$ so is $\dL$, i.e.,
\begin{equation}\label{inv_vol_dL}
L_{\vh}\,\dL=0,
\end{equation}
In the case of a Riemannian surface $(M,g)$ the Liouville volume element on
$h_g^{-1}(1) = \mathcal{S}^{*g}$ is called Riemannian volume element and we denote it by $\dR$.
In this particular case it is easy to check that $\dR$ is invariant by the vertical field $\bs{v}_g$
(actually the Riemannian case toghether with the Lorentzian are the unique ones satisfying the regularity assumptions (\ref{Reg_Assumptions}) for which the canonical vector field $\bs{v}$, defined by relation (\ref{b})
leaves invariant the Liouville volume). Thus, being invariant by $\vh_g$ and $\bs{v}_g$ the Riemannian
volume element is also invariant by their bracket, that is
\begin{equation}\label{inv_vol_dR}
L_{\left[\vh_g , \bs{v}_g \right]}\dR = 0.
\end{equation}
Using relation (\ref{om^d_qom=phidV}), one can easily checked that
for the co-Zermelo problem of the pair $(g,\up)$ the two volume elements $\dL$ and $\dR$ are linked by the relation
\begin{equation}\label{dL=phidR}
{\diffeo{}}^*\dL = \phicoZ\circ \diffeo{-}\, \dR,
\end{equation}
where $\diffeo{}$ is the diffeomorphism defined by relation (\ref{diffeo_F}).

\begin{lemma}\label{GB_for_DZ}
Let $(M,g)$ be a compact, orientable, two-dimensional Riemannian manifold without boundary.
Let $\up$ be a smooth one-form on $M$.
Then,
\begin{equation}\label{GB_neq_rep-DZ} 
\frac{1}{4\pi^2}\int_{\mathcal{S}^{g*}} \kmag^{\sr (g,d\up)} \circ \diffeo{} \,\dR \ge \chi(M),
\end{equation}
\begin{equation}\label{GB_neq_for_DZ} 
\frac{1}{4\pi^2}\int_{\H}\phicoZ \kappa_\coZ^{\sr (g,\up)}\, \dL \ge \chi(M),
\end{equation}
where $\chi(M)$ is the Euler characteristic of the surface $M$.
\end{lemma}
\begin{proof}
According to (\ref{kmag_def}),
\begin{equation*}
\frac{1}{4\pi^2}\int_{\mathcal{S}^{g*}} \kmag^{\sr (g,d\up)} \circ \diffeo{} \,\dR
=   \int_{\mathcal{S}^{g*}} \kappa_g \,\dR
   + \int_{\mathcal{S}^{g*}} \Om^2 \,\dR
   + \int_{\mathcal{S}^{g*}} L_{\left[\vh_g , \bs{v}_g \right]} \Om \,\dR,
\end{equation*}
which, according to the classical Gauss-Bonnet formula and relation (\ref{inv_vol_dR}), is equivlent to
\begin{equation}\label{inter1} 
\frac{1}{4\pi^2}\int_{\mathcal{S}^{g*}} \kmag^{\sr (g,d\up)} \circ \diffeo{} \,\dR
 =  4\pi^2\chi(M) + 2\pi\int_{M} \Om^2 \, dV_g
\ge 4\pi^2\chi(M).
\end{equation}
This proves relation (\ref{GB_neq_rep-DZ})
According to relations (\ref{kmag=phi-rep_coZ}) and (\ref{inv_vol_dL}), we have
\begin{eqnarray}\label{inter2} 
\int_{\HDZ}\phicoZ \mathcal{S}(\phicoZ) \,\dL
&=& - \int_{\HDZ} L_{\vhDZ}\left(\frac{L_{\vhDZ}\phicoZ}{2}\right) \, \dL
    + \int_{\HDZ} \left(\frac{L_{\vhDZ}\phicoZ}{2}\right)^2 \,\frac{\dL}{\phicoZ} \nonumber\\
&=&  \int_{\HDZ} \left(\frac{L_{\vhDZ}\phicoZ}{2}\right)^2 \,\frac{\dL}{\phicoZ}
\ge 0.
\end{eqnarray}
Relation (\ref{GB_neq_for_DZ}) follows from (\ref{inter1}) and (\ref{inter2}), which completes the proof.
\end{proof}

\begin{theorem}\label{GB}
Let $M$ be a compact, orientable, two-dimensional Riemannian manifold without boundary.
If $\kappa$ is the curvature of a Zermelo'-like problem
then,
there exists a canonically defined positive function $\phi$ which is identically equal to one if and only if the problem is Riemannian
such that
\begin{equation}\label{GB_neq} 
\frac{1}{4\pi^2}\int_{\H}\phi \kappa\, \dL \ge \chi(M).
\end{equation}
Moreover, when $\phi$ is identically equal to one relation (\ref{GB_neq}) is the classical Gauss-Bonnet formula.
\end{theorem}
\begin{proof}
It follows straightforwardly from the previous lemma and Proposition \ref{EquivZerDZ}.
\end{proof}\\

It immediately follows from the above theorem that
\begin{theorem}\label{CoroGBInequality}
Zermelo's problems having non positive not identically zero curvature do not exist on two-dimensional tori.
\end{theorem}
\begin{proof}
We prove the result by contradiction. Let $\kappa$ be the curvature of a Zermelo'-like problem on a two-dimensional Riemannian torus and let $\phi$ be the function of Theorem (\ref{GB}). Suppose that $\kappa \le 0$. Since $\kappa$ does not vanish identically, there exists a point $\la \in \H$ such that $\kappa(\la) < 0$, which, in addition with the fact that $\phi$ is a strictly positive function implies that
$\int_{\H}\phi\kappa\,\dL < 0.$
But this contradicts the Gauss-Bonnet inequality of Theorem (\ref{GB}) which, in this case reads
$\int_{\H}\phi\kappa\,\dL \ge 4\pi^2\chi(\T^2) = 0.$
\end{proof}\\
\begin{remark}
{\rm
Although the previous theorem is an immediate consequence of inequality (\ref{GB_neq}), we want to point out that this theorem also follows from a more general fact if ``non-positive" is replaced by ``negative" in its formulation.
Indeed, the flow generated by the Hamiltonian of a smooth control system having negative curvature is Anosov (see \cite{AAAChtch2005}).
Moreover, in the appendix to the paper by Anosov and Sinai \cite{AnosovSinai1967}, Margulis proved that if an Anosov flow operates on a three-dimensional manifold then, its fundamental group has exponential growth.
Therefore, an Anosov flow cannot be carried by a three-dimensional torus since the fundamental group of the latter is the free abelian group
$\mathbb{Z}^3$ which is known to have polynomial and not exponential growth (see e.g. \cite{Gromov1981}).
Finally, one easily check that the hypersurface $\H$ of a Zermelo'-like problem 
(of course, whose drift has Riemannian norm strictly smaller that one)
over a two-dimensional torus
is diffeomorphic to a three-dimensional torus.  
}
\end{remark}

It's not worth mentioning that the Gauss-Bonnet (\ref{GB_neq}) inequality becomes an equality {\sc not only} if the problem is Riemannian.
Indeed,
\begin{proposition}\label{GBequality}
The Gauss-Bonnet inequality of Theorem \ref{GB} is an equality if and only if the drift is identically zero or
the Gaussian curvature of the manifold is zero and the drift has constant Riemannian norm.
\end{proposition}
%
\begin{proof}
It follows from Proposition \ref{EquivZerDZ} that it is enough to prove the result for the co-Zermelo problem of the pair $(g,\up)$.
Let $M = \cup_\alpha O_\alpha$ where the $O_\alpha$'s are domains of local $g$-orthonormal frames
and let $(\bs{e}_1,\bs{e}_2)$ be such a frame.
From relation (\ref{inter1}) we know that
\begin{equation}\label{IntRhoKappa=blablabla}
\int_{\H}\phicoZ \kappa_\coZ^{\sr (g,\up)} \,\dL 
= 4\pi^2\chi(M) + 2\pi\int_{M}\Om^2\, dV_g
+ \int_{\H} \bigg( \frac{L_{\vhmag}\phicoZ}{2\phicoZ} \bigg)^2 {\diffeo{-}}*\dR
\end{equation}
so that the Gauss-Bonnet inequality becomes an equality if and only if
\begin{equation}\label{Om=Lhhatrho=0}
\Om=0 \quad {\rm and} \quad L_{\vhmag}\phicoZ=0
\end{equation}
identically.
On the one hand, the condition $\Om=0$ means that the drift form $\up$ is closed (recall that $\Om$ was defined by $d\up = \Om\, dV_g$),
which implies
\begin{eqnarray}\label{Wclosed}
0 = d\up(\bs{e}_1,\bs{e}_2) &=& L_{\bs{e}_1}\up_2 - L_{\bs{e}_2}\up_1 - \up_1 c_1 - \up_2 c_2,
\end{eqnarray}
where $\up_1 = \lang \up,\bs{e}_1\rang$ and $\up_2 = \lang \up,\bs{e}_2\rang$.

On the other hand, keeping in mind that $\Om=0$ holds true,
condition $L_{\vhmag}\phicoZ = 0$ reads $L_{\diffeo{}_*\vh_g}\phicoZ = 0$.
According to the notations of Example \ref{Riem_geo_pb}, that is
\begin{eqnarray}
0 &=& L_{\cos\theta\bs{e}_1 + \sin\theta\bs{e}_2 + (c_1\cos\theta+c_2\sin\theta)\parfrac{}{\theta}}
(1+\up_1\cos\theta+\up_2\sin\theta) \nonumber \\
&=& ( L_{\bs{e}_1}\up_1 + c_1\up_2 )\cos^2\theta + ( L_{\bs{e}_2}\up_2 - c_2\up_1 )\sin^2\theta \nonumber \\
&&        + ( L_{\bs{e}_1}\up_2 + L_{\bs{e}_2}\up_1 - c_1\up_1 + c_2\up_2 )\cos\theta\sin\theta \label{hriemrho}
\end{eqnarray}
Equations (\ref{Wclosed}) and (\ref{hriemrho}) are thus equivalent to the system of equations
\begin{eqnarray*}
L_{\bs{e}_1}\up_2 - L_{\bs{e}_2}\up_1 - c_1\up_1 - c_2\up_2 &=& 0 \\
L_{\bs{e}_1}\up_1 + c_1\up_2 &=& 0 \\
L_{\bs{e}_2}\up_2 - c_2\up_1 &=& 0 \\
L_{\bs{e}_1}\up_2 + L_{\bs{e}_2}\up_1 - c_1\up_1 + c_2\up_2 &=& 0.
\end{eqnarray*}
Replacing the first and last equations respectively by there sum and difference we equivalently get
\begin{eqnarray}
L_{\bs{e}_1}\up_2 - c_1\up_1 &=& 0 \label{S1} \\
L_{\bs{e}_1}\up_1 + c_1\up_2 &=& 0 \label{S2} \\
L_{\bs{e}_2}\up_2 - c_2\up_1 &=& 0 \label{S3} \\
L_{\bs{e}_2}\up_1 + c_2\up_2 &=& 0.\label{S4}
\end{eqnarray}
Now we differentiate equation (\ref{S4}) along $\bs{e}_1$ and subtract it the differentiation along $\bs{e}_2$
of equation (\ref{S2}). According to (\ref{riemcurv}), we get
\begin{eqnarray}\label{betaK}
0 &=& L_{\bs{e}_1}(\ref{S4}) - L_{\bs{e}_2}(\ref{S2}) \nonumber \\
&=& L_{\bs{e}_1}\circ L_{\bs{e}_2}\up_1 + c_2 L_{\bs{e}_1}\up_2 + \up_2 L_{\bs{e}_1}c_2
       -L_{\bs{e}_2}\circ L_{\bs{e}_1}\up_1 - c_1 L_{\bs{e}_2}\up_2 - \up_2 L_{\bs{e}_2}c_1 \nonumber \\
 &=& L_{[\bs{e}_1,\bs{e}_2]}\up_1 + \up_2(L_{\bs{e}_1}c_2-L_{\bs{e}_2}c_1)
    +(\up_1 c_1)c_2-(\up_1 c_2)c_1  \nonumber \\
 &=& c_1L_{\bs{e}_1}\up_1 + c_2L_{\bs{e}_2}\up_1 + \up_2(L_{\bs{e}_1}c_2-L_{\bs{e}_2}c_1) \nonumber \\
&=& \up_2(-c_1^2-c_2^2+L_{\bs{e}_1}c_2-L_{\bs{e}_2}c_1) \\
&=& \up_2\kappa_g.
\end{eqnarray}
In the same way, using this time equations (\ref{S1}) and (\ref{S3}) we get
\begin{equation}\label{alphaK}
0 = L_{\bs{e}_2}(\ref{S1}) - L_{\bs{e}_1}(\ref{S3}) = \up_1\kappa_g.
\end{equation}

If the Gaussian curvature is identically equal to zero then, the Riemannian manifold is a flat torus. In this case we can chose local coordinates $(q_1,q_2)$ on $M$ such that $\bs{e}_1=\parfrac{}{q_1}$ and $\bs{e}_2=\parfrac{}{q_2}$. In these coordinates equations (\ref{S1}), (\ref{S2}), (\ref{S3}) and (\ref{S4}) read
\begin{equation*}
L_{\bs{e}_i}\up_j = 0, \quad i,j = 1,2,
\end{equation*}
which obviously implies that the coefficients $\up_1$ and $\up_2$ are constant.
Therefore $\up$ has constant Riemannian norm.

If the Gaussian curvature is not identically equal to zero then, it follows from equations (\ref{betaK}) and (\ref{alphaK}) that the form $\up$ must be zero wherever $\kappa_g$ is different from zero.
Consider the set $A = \{ q \in M \, : \, \kappa_g(q)=0 \}$.
If the interior of $A$ is empty it follows from its continuity that $\up$ vanishes identically on M.
If the interior of $A$ is non empty, a similar reasoning as above (done on successively on each domain $O_\alpha$) and the continuity of $\up$ imply that $\up$ has constant Riemannian norm in restriction to the closure of the interior of $A$.
But, $\up|_{M \setminus A} = 0$ and since $(M \setminus A) \cap {\rm clo\,}{\rm int\,}A \neq \emptyset$, by continuity we must have $\up=0$ identically on $M$.
This completes the proof of the theorem.
\end{proof}\\
\begin{remark}\label{Lhphi=0=>Om=0}
{\rm
Equations (\ref{Om=Lhhatrho=0}) are indeed equivalent to the unique equation $L_{\vhmag}\phicoZ = 0$. Namely, $L_{\vhmag}\phicoZ$ is a polynomial of degree two in $\cos\theta$, $\sin\theta$. In particular we have
\begin{eqnarray*}
0 &=& L_{\vhmag}\phicoZ(\pi/2,q) + L_{\vhmag}\phicoZ(-\pi/2,q) = \Om \up_1 \\
0 &=& L_{\vhmag}\phicoZ(0,q) + L_{\vhmag}\phicoZ(\pi,q) = \Om \up_2,
\end{eqnarray*}
which obviously implies that $\Om = 0$.
}
\end{remark}

\section{A E. Hopf theorem for control systems}\label{SectionHopf}
It is well known that Riemannian tori without conjugate points are flat.
This theorem was first proved by E. Hopf in 1943 for the two-dimensional case (see \cite{hopf1948}) and for higher
dimensional manifolds it was proved by D. Burago and S. Ivanov in 1994 (see \cite{burago}). We give in this section a generalization of Hopf's for control systems.

\subsection{Jacobi curves}\label{JacCurvChapGlo}
We introduce here the {\it Jacobi curves}  which are a generalization of the space of Jacobi fields along Riemannian geodesics. Since the construction of Jacobi curves does not depend on the dimension of the manifold, we begin with the general case to then go to our special low-dimensional case.

Let $h$ be the Hamiltonian function of PMP for a time-optimal smooth control problem and $\H$ its hypersurface $h^{-1}(1)$. Let $\e{t}{\vh}:\H \to \H$ denote the flow generated by the Hamiltonian field of PMP $\vh$. This flow defines a one-dimensional foliation $\mathcal{F}$ of $\mathcal{H}$ whose leaves, the trajectories of $\vh$, are transverse to the fibers $T^*_qM$, $q \in M$. This foliation enable us to make the following symplectic reduction.

Consider the canonical projection
\begin{equation*}
\bar{\pi} :\H \to \Sigma=\H/\F.
\end{equation*}
The quotient space $\Sigma$, space of trajectories of $\vh$, is, at least locally, a well-defined smooth manifold and carries a structure of symplectic manifold with symplectic form $\bar{\sigma}$ characterized by the property that its pull-back to $\H$ is the restriction $\sigma|_\H$.

Let $\Pi \subset T\H$ denote the vertical distribution, i.e., $\Pi_\la = T_\la\H_{\pi(\la)}$, $\la \in \H$.
The curve
\begin{eqnarray*}
J_\lambda : \R &\to& T_{\bar{\pi}(\la)}\Sigma \\
t &\mapsto& J_\la(t) = \bar{\pi}_* \circ \e{-t}{\vh}_* \Pi_{\e{t}{\vh}(\lambda)},
\end{eqnarray*}
is called {\it Jacobi curve at $\la$}.
Because the Hamiltonian flow preserves the symplectic structure, it is easy to check that the spaces $J_\la(t)$, $t \in \R$, are Lagrangian subspaces of the
symplectic space $T_{\bar{\pi}(\la)}\Sigma$ so that the Jacobi curves are curves in the Lagrangian Grassmannian
$L(T_{\bar{\pi}(\la)}\Sigma)$.

Recall that the Lagrangian Grassmannian $L(T_{\bar{\pi}(\la)}\Sigma)$ of the symplectic space $T_{\bar{\pi}(\la)}\Sigma$ is defined by:
\begin{equation*}
L(T_{\bar{\pi}(\la)}\Sigma) = \{\Lambda\subset T_{\bar{\pi}(\la)}\Sigma \ |\ \Lambda^\angle=\Lambda\},
\quad \Lambda^\angle = \{\xi \in T_{\bar{\pi}(\la)}\Sigma\ |\ \bar{\sigma}(\xi,\Lambda)=0 \}.
\end{equation*}
The Lagrangian Grassmannian of a symplectic space is a well-defined smooth and compact manifold.
In our particular case of a two-dimensional manifold $M$, the Lagrangian Grassmannian $L(T_{\pi(\la)}\Sigma)$ is diffeomorphic to the
one-dimensional real projective space $\R\mathbb{P}(1)$. Moreover, since the vertical distribution $\Pi$ is generated by the vertical vector field $\bs{v}$ the Jacobi curve can written as
\begin{equation}\label{JacCurv=Rexpadv}
J_\la(t) = \R \left( \bar{\pi}_*\e{t}{\ad \vh}\bs{v}(\la) \right).
\end{equation}

We say that a point $\e{t}{\vh}(\la)$ is {\it conjugate} to $\la$ (or time $t$ is conjugate to zero) if
\begin{equation*}
J_{\la}(t) \cap J_{\la}(0) \neq \{0\}.
\end{equation*}

Most of the material presented in this section can be fund in great details in the papers \cite{AAA_Jac2, AAAGam_Jac1, AAAZel_Jac1}.

\subsection{A E. Hopf theorem}
In this section we prove the following
\begin{theorem}\label{E.Hopf_Th}
Consider a control system $\dot{q}=\bs{f}(q,u)$ on a compact surface $M$ without boundary. Assume that the curves of admissible velocities are strongly convex curves surrounding the origin.
Then, if there is no conjugate points on $M$ the total curvature
$\int_{\mathcal{H}}\kappa\, \dL$ must be negative or zero. In the latter case $\kappa$ must be zero.
\end{theorem}
\begin{proof}
Notice that because the curves of admissible velocities are strongly convex curves surrounding the origin, the manifold $\H$ is compact.
Although the proof we make here essentially follows the one given by Hopf in \cite{hopf1948}, it will however be exposed in a more intrinsic and geometrical manner. The first step in the proof consists in the construction of a well-defined function on  any extremal of our system, i.e., a function that does not depend on time but only on the point of the extremal. To do so we use the notion of Jacobi curve described in the previous section.

Let $\la$ be a point of the hypersurface $\H \subset T^*M$ and let $J_\la(t)$ be the Jacobi curve associated with the extremal $\e{t}{\vh}(\la)$.
we have
\begin{equation*}
J_\la(t) = \R \left( \bar{\pi}_*\e{t}{\ad \vh}\bs{v}(\la) \right) \in \R\mathbb{P}(1),
\end{equation*}
with
\begin{equation*}
\e{t}{\ad \vh}\bs{v}(\la) = \beta(t,\la)\bs{v}(\la) + \gamma(t,\la) \left[ \bs{v} , \vh \right](\la).
\end{equation*}
Considering $(\beta:\gamma)$ as homogeneous coordinate in $\R\mathbb{P}(1)$, we can identify the Jacobi curve with the curve
\begin{equation*}
t \mapsto (\beta(t,\la):\gamma(t,\la)).
\end{equation*}
From the non existence of conjugate points it follows that $\gamma(t,\la) \neq 0$ for $t \neq 0$. We can thus use the chart
$(\beta:\gamma)\mapsto \frac{\beta}{\gamma}$ and make the identification
\begin{equation*}
J_\la(t) = y_t(\la) = \frac{\beta(t,\la)}{\gamma(t,\la)}, \quad t \neq 0.
\end{equation*}
It turns out (see e.g. \cite{AAAbook, ulysseThese}) that the coefficients $\beta$ and $\gamma$ are solutions of the Cauchy problems
\begin{eqnarray*}
&& \ddot{\beta\,}\!+\kappa_t\beta = 0, \quad \beta(0)=1, \quad \dot{\beta\,}\!(0)=0, \quad \kappa_t=\kappa(\e{t}{\vh}(\la)),\\
&& \ddot{\gamma\,}\! +\kappa_t\gamma =0, \quad \gamma(0)=0, \quad \dot{\gamma\,}\!(0)=1,
\end{eqnarray*}
which shows in particular that $\beta$ and $\gamma$ are two linearly independent solutions of the Hill equation $\ddot{x}+\kappa_t x=0$.
The derivative with respect to time of the function $y_t$ is
\begin{equation*}
\frac{d y_t}{dt} = \frac{\dot{\beta\,}\!\gamma - \beta \dot{\gamma\,}\!}{\gamma^2}
\end{equation*}
and because the Wronskian
\begin{equation*}
\dot{\beta\,}\!(0,\la)\gamma(0,\la) - \beta(0,\la) \dot{\gamma\,}\!(0,\la)= -1,
\end{equation*}
the function $y_t$ is strictly decreasing or, equivalently the Jacobi curve is strictly decreasing in $\R\mathbb{P}(1)$. Since $y_t$ is strictly decreasing
its limit as $t$ goes to infinity exists. Moreover, because of the non existence of conjugate points, this limit is finite.
Indeed, notice that because of the initial conditions $\beta(0,\la)=1$, $\gamma(0,\la)=0$ and $\dot{\gamma\,}\!(0,\la)=1$ we have for $t$ small enough
\begin{equation}\label{utsmall}
y_t(\la)>0, \quad y_{-t}(\la)<0.
\end{equation}
So if we suppose that
\begin{equation}\label{limu=infinity}
\lim_{t \to +\infty}y_t(\la) = -\infty,
\end{equation}
it would follow from Equations (\ref{utsmall}) and from the strict monotonicity of $y_t$
the existence of $t^{\sr -} < 0 < t^{\sr +}$ such that $y_{t^{\sr -}}(\la)=y_{t^{\sr +}}(\la)$. Then, the time reparametrization $\tau=t-t^{\sr -}$ would imply that time $\tau=t^{\sr +} - t^{\sr -}$ is conjugate to $\tau = 0$, which is a contradiction. Hence, the function $y^{\sr +}$ defined by
\begin{equation*}
y^{\sr +}(\la) = \lim_{t \to +\infty} y_t(\la), \quad \la \in \H
\end{equation*}
is a well defined function on the manifold $\H$.
Equivalently, the distribution $\Pi_{\la}^{\sr +} \in T\H$ defined by
\begin{equation*}
\Pi_{\la}^{\sr +} = \lim_{t \to +\infty} J_{\la}(t) = \R \Big( y^{\sr +}\bs{v} + \left[ \bs{v} , \vh \right] \Big)
\end{equation*}
is a well defined distribution on $\H$ transverse to the vertical distribution.
This distribution $\Pi_{\la}^\infty$ is, by definition, invariant by the flow of $\vh$.
In terms of function $y^{\sr +}$,
this invariance reads
\begin{equation*}
\Big[ \vh, y^{\sr +}\bs{v} + \left[ \bs{v} , \vh \right] \Big] = \alpha \left( y^{\sr +}\bs{v} + \left[ \bs{v} , \vh \right] \right),
\end{equation*}
or, equivalently
\begin{equation}\label{huv}
L_{\vh}y^{\sr +}\bs{v} + y^{\sr +}\Big[ \bs{v} , \vh \Big] + \Big[ \vh , \Big[ \bs{v} , \vh \Big] \Big]
= \alpha y^{\sr +}\bs{v} + \alpha \Big[ \bs{v} , \vh \Big],
\end{equation}
where $\alpha$ is function on $\mathcal{H}$. Solving (\ref{huv}) for $\alpha$ gives
\begin{equation*}
\alpha=-y^{\sr +}\quad {\rm and}\quad L_{\vh}y^{\sr +} + \kappa-\alpha y^{\sr +}=0,
\end{equation*}
which shows that $y^{\sr +}$ satisfies the Riccati equation
\begin{equation}\label{Ric}
L_{\vh}y^{\sr +} + y^{{\sr +}2}+\kappa=0.
\end{equation}
As a limit of smooth functions, $y^{\sr +}$ is clearly measurable. $y^{\sr +}$ is also uniformly bounded as shows lemma 2.1 of \cite{Green1954} and thus it is integrable.
If we now integrate equation (\ref{Ric}) over $\H$ with respect to the Liouville volume $\dL$, the first term in the left-hand side of (\ref{Ric}) will disappear since the Liouville volume is invariant by the flow of $\vh$. As a result we obtain
\begin{equation}\label{totalK}
\int_\mathcal{H}\kappa\, \dL=-\int_\mathcal{H}y^{{\sr +}2}\, \dL
\end{equation}
which immediately proves the validity of the first part of the theorem. If we now suppose that the total curvature $\int_\H \kappa\, \dL$ is zero it follows from (\ref{totalK}) that the function $y^{\sr +}$ must vanish everywhere on $\H$. According to (\ref{Ric}) $\kappa$ must therefore vanish everywhere.
\end{proof}\\

We say that a control system $\dot{q}=\bs{f}(q,u)$ {\it is flat} if it is feedback equivalent to a control system of
the form $\dot{q}=\bs{f}(u)$.

In the Riemannian case, a direct consequence of the Gauss-Bonnet and Theorem \ref{E.Hopf_Th} is that two-dimensional Riemannian tori without conjugate points are flat.
Contrary to the Riemannian situation, we shall see that Zermelo'-like problems without conjugate points on tori are not necessarily flat.

The following three theorems give us a well understanding of the Zermelo'-like situation.
To simplify notations, we omit the pair $(g,\up)$ in the writing of curvature and, the diffeomorphism (\ref{diffeo_F}) in formulas since, anyway, its action is clear.
%
\begin{theorem}\label{E.Hopf_Th_For_Zer}
Consider a co-Zermelo problem on a compact Riemannian surface without boundary.
If there is no conjugate points then the total curvatures
$\int_{\HDZ}\kappa_\coZ \dL$
and
$\int_{\mathcal{S}^{g*}} \kmag \dR$
have to be negative or zero.
In the latter case the considered co-Zermelo problem is flat.
\end{theorem}
\begin{proof}
The part of the theorem concerning $\kappa_\coZ$ is given by Theorem \ref{E.Hopf_Th}.
In order to check that the same conclusion holds for the curvature $\kmag$, let us see how changes the function $y^{\sr +}$ constructed in the proof of Theorem \ref{E.Hopf_Th} under a reparametrization. For simplicity, denote $\psi^2 = \varphi$.
In a general manner, we have
\begin{equation*}
\vh=\frac{\he}{\psi^2}\quad {\rm and}\quad \bs{v}=\psi\ve,
\end{equation*}
and we compute the new function $\hat{y}^{\sr +}$:
\begin{eqnarray*}
y^{\sr +}\bs{v} + \Big[ \bs{v},\vh \Big]
&=& y^{\sr +}\psi\ve + \Big[ \psi\ve + \frac{1}{\psi^2}\he \Big]
 =  y^{\sr +}\psi\ve + \frac{1}{\psi}\Big[\ve,\he\Big]
   - \frac{1}{\psi^2}L_{\he}\psi\ve \quad \left( {\rm mod\ }\vh \right) \nonumber \\
&=& \left( y^{\sr +}\psi-L_{\vh}\psi\right)\ve + \frac{1}{\psi}\Big[\ve,\he\Big] \quad \left( {\rm mod\ }\vh \right).
\end{eqnarray*}
We thus have
\begin{equation}\label{y+hat}
\Pi^{\sr +} = \R\Big(\hat{y}^{\sr +}\ve+\Big[\ve,\he\Big]\Big),
\quad \hat{y}^{\sr +} = y^{\sr +}\psi^2 - \psi L_{\vh}\psi.
\end{equation}
In the same way as for the function ${y}^{\sr +}$ it is easy to see that the function $\hat{y}^{\sr +}$ satisfies the Riccati equation
\begin{equation}\label{Ric_for_y+hat}
L_{\he}\hat{y}^{\sr +}+\hat{y}^{{\sr +}2}+\hat{\kappa} = 0.
\end{equation}
Notice that the Riemannian volume element $\dR$ is invariant by $\vhmag$ since
\begin{equation*}
L_{\vhmag} \dR = L_{\Om \bs{v}_g} \dR = d(\Om dV) = d(d\up) = 0.
\end{equation*}
Therefore the integration of (\ref{Ric_for_y+hat}) leads to
\begin{equation*}
\int_{\mathcal{S}^{g*}} \kmag \dR = - \int_{\mathcal{S}^{g*}} \hat{y}^{{\sr +}2} \dR \le 0.
\end{equation*}
This prove the first part of the theorem and a similar argument as the one used in the proof of Theorem \ref{E.Hopf_Th} shows that $\kmag$ is zero everywhere when $\int_{\mathcal{S}^{g*}} \kmag \,\dR = 0$.

We now complete the proof showing that the co-Zermelo problem is flat when the total curvatures
$\int_{\HDZ} \kappa_\coZ \,\dL$
and
$\int_{\mathcal{S}^{g*}} \kmag \,\dR$
are both zero.
In that case, we must have $\kmag = 0$ and $\kappa_\coZ = 0$ everywhere.
In particular it implies
\begin{equation*}
0 = \int_{\mathcal{S}^{g*}} \kmag \,\dR = \int_{\H} \varphi\kappa_\coZ\, \dL =0,
\end{equation*}
i.e. (see the proof of Lemma \ref{GB_for_DZ}),
\begin{equation*}
0 = 4\pi^2 \chi(M) + 2\pi\int_M\Om^2dV_g
  = 4\pi^2 \chi(M) + 2\pi\int_M\Om^2dV_g + \int_{\H}\left(\frac{L_{\vhDZ}\varphi}{2}\right)^2 \, \frac{\dL}{\varphi},
\end{equation*}
which is equivalent to
$L_{\vhDZ}\varphi = 0$.
Therefore, according to Theorem \ref{GBequality} and Remark \ref{Lhphi=0=>Om=0},
or the form $\up$ is different from zero and in this case the conclusion is obtained,
or the form $\up$ is identically zero and in this case the problem is Riemannian.
In the latter case, we have $0 = \kmag = \kappa_g$ which,
on the one hand, implies that the Riemannian surface is flat and,
on the other hand, according to the Gauss-Bonnet formula, it implies that the surface is a torus.
The proof is complete.
\end{proof}

A direct consequence of this theorem are the following corollaries.
\begin{coro}\label{E.Hopf_Th_for_DZ_on_T^2}
If a co-Zermelo problem on a two-dimensional Riemannian torus has no conjugate points
then, the torus is flat and the drift one-form is closed.
In particular, time-optimal trajectories are straight lines. 
\end{coro}
%
\begin{coro}\label{Flat_Zer'-like_pbs_on_T^2}
Zermelo'-like problems without conjugate points on two-dimensional Riemannian tori are flat
if and only if their total curvature is zero. 
\end{coro}
%

\subsection{A natural question}

In the proof of Theorem \ref{E.Hopf_Th} we constructed a function $y^{\sr +}$ well-defined on $\H$ that satisfies Riccati equation (\ref{Ric}). This construction is valid along every regular extremal without conjugate points. Recall moreover that a control system with negative curvature does not admit conjugate points. A very natural question is thus the following:\\

\noindent {\it Does a control system without conjugate points admits a non positive $\varphi$-reparametrized curvature?}\\

Since the function $\hat{y}^{\sr +}$ satisfies Riccati equation (\ref{Ric_for_y+hat}), the question can be reformulated
in the following manner: {\it does there exists a non vanishing function $\psi$, say $\psi > 0$ for simplicity, such that $L_{\he}\hat{y}^{\sr +}=0$, or equivalently such that $L_{\vh}\hat{y}^{\sr +}=0$?}
According to relation (\ref{y+hat}),
\begin{eqnarray*}
L_{\vh}\hat{y}^{\sr +}
&=& L_{\vh}(y^{\sr +}\psi^2 - \psi L_{\vh}\psi)
= \psi^2 L_{\vh}y^{\sr +} + 2y^{\sr +}\psi L_{\vh}\psi -(L_{\vh}\psi)^2 - \psi L_{\vh}^2\psi,
\end{eqnarray*}
so that (dividing by $\psi^2$) $L_{\vh}\hat{y}^{\sr +}=0$ is equivalent to
\begin{equation*}
L_{\vh}y^{\sr +} + 2y^{\sr +}\left(\frac{L_{\vh}\psi}{\psi}\right)
- \left(\frac{L_{\vh}\psi}{\psi}\right)^2 - \frac{L_{\vh}^2 \psi}{\psi} =0,
\end{equation*}
i.e., to
\begin{equation}\label{R1}
L_{\vh}y^{\sr +} + 2y^{\sr +} L_{\vh}\log \psi - (L_{\vh}\log \psi)^2 - \frac{L_{\vh}^2 \psi}{\psi} =0.
\end{equation}
Denote $g = \log \psi$. We have
\begin{equation*}
L_{\vh}^2g
= L_{\vh}(L_{\vh}\log \psi)
= L_{\vh} \left( \frac{L_{\vh}\psi}{\psi} \right)
= \frac{(L_{\vh}^2 \psi)\psi - (L_{\vh}\psi)^2}{\psi^2}
= \frac{L_{\vh}^2 \psi}{\psi} - (L_{\vh}g)^2,
\end{equation*}
or equivalently
\begin{equation*}
\frac{L_{\vh}^2 \psi}{\psi} = L_{\vh}^2g+(L_{\vh}g)^2.
\end{equation*}
This implies that equation (\ref{R1}) is equivalent to
\begin{equation*}
L_{\vh}y^{\sr +} +2y^{\sr +}L_{\vh}g - 2(L_{\vh}g)^2 - L_{\vh}^2 g=0,
\end{equation*}
i.e., to the Riccati equation
\begin{equation}\label{Ricc}
L_{\vh}z + 2z^2 -2y^{\sr +}z - L_{\vh}y^{\sr +}=0,
\end{equation}
where we have set $z=L_{\vh}g$.

The function $z=y^{\sr +}$ is solution to Riccati equation (\ref{Ricc}). Thus we will have the required reparametrization
of $\vh$ if we can solve the equation
\begin{equation}\label{L2hlogf=u}
L_{\vh}^2\log \psi = y^{\sr +}
\end{equation}
globally on the three-dimensional manifold $\H$.
The first thing we need for the resolution of equation (\ref{L2hlogf=u}) is the continuity of the function $y^{\sr +}$ on $\H$. In the case of hyperbolic systems (see \cite{HasselblattKatok1995} for the definition), the function $y^{\sr +}$ is easily seen to be continuous due to some ``exponential estimates" along the stable distribution (see \cite{HasselblattKatok1995}). Also, for such systems the function $y^{\sr +}$ is in general never differentiable and even never Lipschitz continuous but only H\"older continuous (see \cite{HasselblattKatok1995} Theorem 19.1.6 of Chapter 19). In the case of systems without conjugate points the situation is quite different because we do not have the exponential estimates and by consequence the continuity of the function $y^{\sr +}$ is not so obvious.
What we can ensure is the following.
\begin{lemma}\label{uIsContinuous}
The function $y^{\sr +}$ defined above is upper semi-continuous.
\end{lemma}
\begin{proof}
Let $(\la_n)_{n \in \N} \subset \H$ be a converging sequence to $\la \in \H$.
Since $y_t(\la_n)$ is decreasing in $t$, it follows that
\begin{equation*}
y_t(\la_n) \ge y^{\sr +}(\la_n)= \lim_{t \to +\infty}y_t(\la_n).
\end{equation*}
Taking the $\lim\inf$ as $n$ tends to $+\infty$ in the previous relation, we get since $y_t(\la)$ is continuous in $(t,\la)$
\begin{equation*}
y_t(\la) \ge \liminf_{\la_n \to \la } y^{\sr +},
\end{equation*}
and then, letting $t$ going to $+\infty$ leads to
\begin{equation*}
y^{\sr +}(\la) \ge  \liminf_{\la_n \to \la } y^{\sr +},
\end{equation*}
which proves the upper semi-continuity of $y^{\sr +}$.
\end{proof}\\

Suppose that the function $y^{\sr +}$ is continuous. It implies that we can solve locally equation (\ref{L2hlogf=u}). In order to solve this equation globally, the question is more delicate because the problem is closely related to the fact that the quotient manifold $\Sigma$ (defined in Section \ref{JacCurvChapGlo}) is globally defined.
It is not our scope to discuss this problem here. However we can say the following. Let $\tilde{M}$ be the universal covering of $M$. because of the non existence of conjugate points, $\tilde{M}$ is diffeomorphic to $\R^2$. Let
\begin{equation}\label{LiftedSystem}
\dot{\tilde{q}} = \bs{\tilde{f}}(\tilde{q},u), \quad \tilde{q} \in \tilde{M},\quad u \in U,
\end{equation}
be the lift on $\tilde{M}$ of the control system $\dot{q}=\bs{f}(q,u)$, and $\tilde{\H}$ be the corresponding Hamiltonian hypersurface. Then, the continuity of $y^{\sr +}$ implies that when the control system $\dot{q}=\bs{f}(q,u)$ has no conjugate points then, there exists a reparametrization of $\vh$ or, equivalently globally defined function $\psi$ satisfying equation (\ref{L2hlogf=u}), such that the lifted system (\ref{LiftedSystem}) has negative curvature.
Unfortunately, $y^{\sr +}$ is in general not a continuous function as shown by Ballmann, Brin and Burns in \cite{BBB1987} where they gave an example of a two-dimensional compact surface without conjugate points where $y^{\sr +}$ fails to be continuous.

A nice work would be to characterize two-dimensional smooth control systems without conjugate points where this function fails to be continuous. What is the geometrical property that prevents $y^{\sr +}$ from being continuous?

\section{Conclusion}
We conclude this paper with a brief discussion of the extension of our results to more general structures than Riemannian surfaces.
Of course, Zermelo'-like problems can be defined on any manifold equipped with a geometric structure defined by an optimal control problem of type (\ref{optpb_dyn})-(\ref{optpb_cost}).
A natural class of geometric structures on which generalize our result is the class of manifolds equipped with a Finsler metric (see the book of Chern and Shen \cite{ChernShen} for a nice and brief presentation of Riemann-Finsler geometry).
Unfortunately, since the Gauss-Bonnet formula is not true for any Finsler surfaces results from Section \ref{SectionGB} can not be extented to all of these structures.
One has to limits itself to Zermelo'-like problems on Landsberg surfaces on which almost all results from Section \ref{SectionGB} remain true.
Roughly speaking, a Landsberg surface is a Finsler surface on which the Gauss-Bonnet formula remains true.
(up to change the classical $2\pi$ factor in the formula by the centro-affine length $\ell$ of the curve $\H_q$
which, let us recall it, is defined by $\ell = \int_{\H_q}\mu|_{\H_q}$
where $\mu$ is a one-form on the hypersuface $\H$ such that
$\lang \mu , \bs{v} \rang = 1$
). Without entering into details one can see that the Gauss-Bonnet formula still holds on Landsberg surfaces due to the fact that the centro-affine length of the curve $\H_q$
does not depend on the base point $q$.
This property is characterized by the fact that the invariant $b$ that appears in relation (\ref{b}) is a first integral of the vector field $\vh$ (see \cite{ulysseThese} for details).
If we now consider Zermelo'-like problems on Landsberg surfaces, on the one hand, the Gauss-Bonnet inequality (\ref{GB_neq}) still holds true. The proof is the same but this time one has to be more carefull because the Landsberg volume element $\dL_{\rm land}$ is not invariant under the vertical Landsberg field $\bs{v}_{\rm land}$.
Indeed, one can easily check that $L_{\bs{v}_{\rm land}}\,\dL_{\rm land} = b\,\dL_{\rm land}$.
Anyway, $\dL_{\rm land}$ is still invariant under the bracket $[\vh_{\rm land} , \bs{v}_{\rm land}]$ since
\begin{equation*}
L_{\big[ \vh_{\rm land} , \bs{v}_{\rm land} \big] } \dL_{\rm land} = \big( L_{\vh_{\rm land}}b \big) \dL_{\rm land} = 0.
\end{equation*} 
On the other hand, Theorem \ref{E.Hopf_Th_For_Zer} and its Corollaries \ref{E.Hopf_Th_for_DZ_on_T^2} and \ref{Flat_Zer'-like_pbs_on_T^2} do not generalize to Zermelo'-like problems on Landsberg surfaces.
The reason is the following:
Landsberg surfaces of zero curvature are not necessary flat (see \cite{ulysseThese}, Theorem 4.3.3).

\subsection*{Acknowledgments}
I am grateful to Professor Andrei A. Agrachev for fruitful discussions.
%

\end{document}